\documentclass[12pt]{article}
\usepackage{amsfonts}

\usepackage{mathrsfs}
\usepackage{amscd}
\usepackage{amsmath,amsfonts,amssymb,amscd}
\usepackage{indentfirst,graphics,epsfig,psfrag}
\input{epsf}

\usepackage{amsmath,amssymb,amsthm, amscd, amsfonts, graphicx,ifpdf}
\usepackage{lineno}

\newtheorem{thm}{Theorem}[section]

\newtheorem{lem}{Lemma}[section]

\newtheorem{rema}[thm]{Remark}
\newtheorem{obv}[thm]{Observation}
\def\qed{\nopagebreak\hfill{\rule{4pt}{7pt}}
\medbreak}

\def\pf{\noindent {\it Proof.} }

\title{\bf Sharp bounds for the generalized connectivity $\kappa_3(G)$
\footnote{Supported by NSFC, PCSIRT and the ``973" program.}}

\author{
\small Shasha Li, Xueliang Li, Wenli Zhou\\
\small Center for Combinatorics and LPMC-TJKLC\\
\small Nankai University, Tianjin 300071, China.\\
\small  Email: lss@cfc.nankai.edu.cn,\\
\small lxl@nankai.edu.cn, louis@cfc.nankai.edu.cn\\
}
\date{}
\begin{document}

\maketitle

\begin{abstract}
Let $G$ be a nontrivial connected graph of order $n$ and let $k$ be
an integer with $2\leq k\leq n$. For a set $S$ of $k$ vertices of
$G$, let $\kappa (S)$ denote the maximum number $\ell$ of
edge-disjoint trees $T_1,T_2,\ldots,T_\ell$ in $G$ such that
$V(T_i)\cap V(T_j)=S$ for every pair $i,j$ of distinct integers with
$1\leq i,j\leq \ell$. A collection $\{T_1,T_2,\ldots,T_\ell\}$ of
trees in $G$ with this property is called an internally disjoint set
of trees connecting $S$. Chartrand et al. generalized the concept of
connectivity as follows: The $k$-$connectivity$, denoted by
$\kappa_k(G)$, of $G$ is defined by
$\kappa_k(G)=$min$\{\kappa(S)\}$, where the minimum is taken over
all $k$-subsets $S$ of $V(G)$. Thus $\kappa_2(G)=\kappa(G)$, where
$\kappa(G)$ is the connectivity of $G$.

In general, the investigation of $\kappa_k(G)$ is very difficult. We
therefore focus on the investigation on $\kappa_3(G)$ in this paper.
We study the relation between the connectivity and the
$3$-connectivity of a graph. First we give sharp upper and lower
bounds of $\kappa_3(G)$ for general graphs $G$, and construct two
kinds of graphs which attain the upper and lower bound,
respectively. We then show that if $G$ is a connected planar graph,
then $\kappa(G)-1 \leq \kappa_3(G)\leq \kappa(G)$, and give some
classes of graphs which attain the bounds. In the end we show that
the problem whether $\kappa(G)=\kappa_3(G)$ for a planar graph $G$
can be solved in polynomial time.\\[3mm]
{\bf Keywords:} connectivity, $k$-connectivity, internally disjoint
trees (paths)\\[3mm]
{\bf AMS Subject Classification 2000:} 05C40, 05C05, 05C38
\end{abstract}

\section{Introduction}

We follow the terminology and notations of \cite{Bondy} and all
graphs considered here are always simple. As usual, the union of two
graphs $G$ and $H$ is the graph, denoted by $G\cup H$, with vertex
set $V(G)\cup V(H)$ and edge set $E(G)\cup E(H)$. Let $T$ be a set
of vertices. Then, $G-T$ is the graph obtained from $G$ by deleting
all the vertices in $V(G)\cap T$ together with their incident edges.
A path $P=x_0x_1\ldots x_k$ is called an {\it $x_0x_k$-path},
denoted by $x_0Px_k$. For the $x_0x_k$-path $P$, we denote three
special subpaths of $P$ by $\hat{x}_0P\hat{x}_k:=x_1\ldots x_{k-1}$,
$\hat{x}_0Px_k:=x_1\ldots x_k$ and $x_0P\hat{x}_k:=x_0\ldots
x_{k-1}$. For $X=\{x_1,x_2,\ldots,x_k\}$ and
$Y=\{y_1,y_2,\ldots,y_k\}$, an {\it $XY$-linkage} is defined as a
set of $k$ disjoint paths $x_iP_iy_i$, $1\leq i\leq k$. The
$connectivity$ $\kappa(G)$ of a graph $G$ is defined as the minimum
cardinality of a set $Q$ of vertices of $G$ such that $G-Q$ is
disconnected or trivial. A well-known theorem of Whitney
\cite{Whitney} provides an equivalent definition of connectivity.
For each 2-subset $S=\{u,v\}$ of vertices of $G$, let $\kappa(S)$
denote the maximum number of internally disjoint $uv$- paths in $G$.
Then $\kappa(G)=$min$\{\kappa(S)\}$, where the minimum is taken over
all 2-subsets $S$ of $V(G)$.

In \cite{Chartrand}, the authors generalized the concept of
connectivity. Let $G$ be a nontrivial connected graph of order $n$
and let $k$ be an integer with $2\leq k\leq n$. For a set $S$ of $k$
vertices of $G$, let $\kappa (S)$ denote the maximum number $\ell$
of edge-disjoint trees $T_1,T_2,\ldots,T_\ell$ in $G$ such that
$V(T_i)\cap V(T_j)=S$ for every pair $i,j$ of distinct integers with
$1\leq i,j\leq \ell$. A collection $\{T_1,T_2,\ldots,T_\ell \}$ of
trees in $G$ with this property is called an {\it internally
disjoint set of trees connecting $S$}. The $k$-$connectivity$,
denoted by $\kappa_k(G)$, of $G$ is then defined by
$\kappa_k(G)=$min$\{\kappa(S)\}$, where the minimum is taken over
all $k$-subsets $S$ of $V(G)$. Thus, $\kappa_2(G)=\kappa(G)$.

Chartrand et al. in \cite{Chartrand} proved that if $G$ is the
complete 3-partite graph $K_{3,4,5}$, then $\kappa_3(G)=6$. They
also gave an general result for the complete graph $K_n$:
\begin{thm}\label{thm0}
For every two integers $n$ and $k$ with $2\leq k\leq n$,
\begin{center}
$\kappa_k(K_n)=n-\lceil k/2\rceil$.
\end{center}
\end{thm}

In general, the investigation of $\kappa_k(G)$ is very difficult.
Therefore, in this paper we will focus on the investigation of
$\kappa_3(G)$. We study the relation between the connectivity and
the $3$-connectivity of a graph. First, we give sharp upper and
lower bounds of $\kappa_3(G)$ for general graphs $G$, and construct
two kinds of graphs which attain the upper and lower bound,
respectively. Then, we study the 3-connectivity for the planar
graphs. We will show that if $G$ is a connected planar graph, then
$\kappa(G)-1 \leq \kappa_3(G)\leq \kappa(G)$, and give some classes
of graphs which attain the bounds. In the end, we show that the
problem whether $\kappa(G)=\kappa_3(G)$ for a planar graph $G$ can
be solved in polynomial time.

\section{Upper and lower bounds}

Before we give the main results, there is an easy observation:

\begin{obv}\label{obv1}
If $G'$ is a spanning subgraph of $G$, then $\kappa_{k}(G')\leq
\kappa_{k}(G)$ for $2\leq k\leq n$.
\end{obv}

Now we give an upper bound of $\kappa_3(G)$.

\begin{thm}\label{thm1}
Let $G$ be a connected graph with $n$ vertices. Then
$\kappa_3(G)\leq \kappa(G)$, and moreover, the upper bound is sharp.
\end{thm}
\pf We prove the theorem by three cases on $\kappa(G)$.

{\it Case 1:} $\kappa(G)=n-1$.

Then $G$ must be a complete graph $K_n$. By Theorem \ref{thm0}, we
know $\kappa_3(K_n)=n-\lceil \frac{3}{2}\rceil=n-2$. So
$\kappa_3(G)=n-2\leq \kappa(G)=n-1$.

{\it Case 2:} $\kappa(G)=n-2$.

Let $Q$ be an $(n-2)$-vertex cut of $G$. Here and in what follows,
by a {\it $k$-vertex cut} we mean a vertex cut that have $k$
vertices. Assume $V(G)-Q=\{u,v\}$ such that $u$ and $v$ are two
nonadjacent vertices and both of them are adjacent to all vertices
in $Q$. If $Q$ is a clique, it is easy to check that
$\kappa_3(G)=n-2$. Otherwise, $G$ must have a spanning supergraph
$G'=K_n-uv$ (i.e., $G$ is a spanning subgraph of $G'$). By
Observation \ref{obv1}, we get $\kappa_3(G)\leq
\kappa_3(G')=n-2=\kappa(G)$.

{\it Case 3:} $1\leq \kappa(G)\leq n-3$.

Let $Q$ be a $\kappa(G)$-vertex cut of $G$. Then $G-Q$ has at least
2 components. Since $|Q|\leq n-3$, we can choose a vertex set $S$
consisting of three vertices which are not in $Q$, such that two of
the three vertices are in different components. Then we know that
any tree connecting $S$ must contain a vertex in $Q$. By the
definition of $\kappa(S)$, we get $\kappa_3(G)\leq \kappa(S)\leq
|Q|=\kappa(G)$.

From the above, we conclude that $\kappa_3(G)\leq \kappa(G)$.

Furthermore, for any two integers $k\geq 1$ and $n\geq k+2$,
consider the graph $G=K_k \bigvee(n-k)K_1$. Then, obviously
$\kappa(G)=k$, and it is not difficult to check that
$\kappa_3(G)=k$. So $\kappa_3(G)=\kappa(G)=k$, and therefore the
upper bound is sharp.\qed

In the following, we will give a lower bound of $\kappa_3(G)$.
Before proceeding, we recall the {\it Fan Lemma}, which will be used
frequently in the sequel.

\begin{lem} (The Fan Lemma \cite{Bondy}) \label{lem1}
Let $G$ be a $k$-connected graph, $x$ a vertex of $G$, and let
$Y\subseteq V-\{x\}$ be a set of at least $k$ vertices of $G$. Then
there exists a $k$-fan in $G$ from $x$ to $Y$, namely there exists a
family of $k$ internally disjoint $(x,Y )$-paths whose terminal
vertices are distinct in $Y$.
\end{lem}

Our lower bound is given as follows:

\begin{thm}\label{thm2}
Let $G$ be a connected graph with $n$ vertices. For every two
integers $k$ and $r$ with $k\geq0$ and $r\in\{0,1,2,3\}$, if
$\kappa(G)=4k+r$, then $\kappa_3(G)\geq 3k+\lceil
\frac{r}{2}\rceil$. Moreover, the lower bound is sharp.
\end{thm}

Before proving the theorem, we need some preparations. Denote
$\kappa(G)$ by $\kappa$ for short. First, we introduce an operation
called ``Path-Transformation", which can adjust paths in order to
attain some structure we want. More explicitly, first we are given
$\kappa$ $v_1v_2$-paths $P_1,P_2,\ldots,P_{\kappa}$ such that $v_3$
is on $t$ paths $P_1,\ldots,P_t$ of them for some $1\leq t < \lceil
\frac{\kappa}{2}\rceil$, and except $v_3$ the $\kappa$ paths have no
internal vertices in common. For $X=V(P_{t+1}\cup \cdots\cup
P_{\kappa})$, by a family of $\kappa$ internally disjoint
$(v_3,X)$-paths and the ``Path-Transformation" , we adjust the paths
$P_1,\ldots,P_t$ to get $\kappa$ $v_1v_2$-paths
$P_1',\ldots,P_{t}',P_{t+1},\ldots,P_{\kappa}$ which still have the
former structure, and in addition, there is a family of $\kappa-2t$
internally disjoint $(v_3,X)$-paths avoiding the vertices in
$V(P_1'\cup \cdots\cup P_{t}'-\{v_1,v_2,v_3\})$. The following
Figure \ref{fig1} shows the ``Path-Transformation''.

\begin{figure}[h,t]
\begin{center}
\input{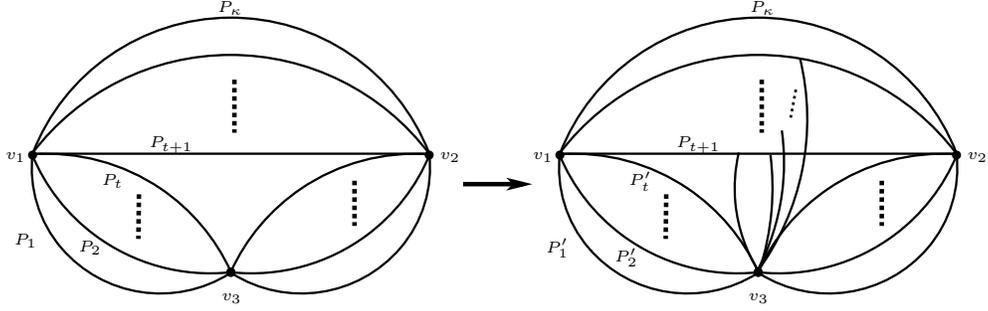}
\caption{Illustration for ``Path Transformation"} \label{fig1}
\end{center}
\end{figure}

Now, we mainly describe how to adjust paths and why the operation
can get the structure we want.

Let $P_1,P_2,\ldots,P_{\kappa}$ be $\kappa$ $v_1v_2$-paths such that
$v_3$ is on $t$ paths $P_1,\ldots,P_t$ of them for some $1\leq t <
\lceil \frac{\kappa}{2}\rceil$, and the $\kappa$ paths have no
internal vertices in common except $v_3$. Then let $X=V(P_{t+1}\cup
\cdots\cup P_{\kappa})$. Since $G$ is $\kappa$-connected and if
$|X|\geq \kappa$ (the case $|X|< \kappa$ will be illustrated later),
there is a $\kappa$-fan $\{M_1,M_2,\ldots,M_{\kappa}\}$ from $v_3$
to $X$. Let $V(M_1\cup \cdots\cup M_{\kappa})\cap V(P_1\cup
\cdots\cup P_{t}-\{v_1,v_2\})=N$, where $N\neq \emptyset$, since at
least the vertex $v_3$ belongs to $N$. $P_1,\ldots,P_t$ can be
regarded as $2t$ paths
$v_1P_1v_3,v_1P_2v_3,\ldots,v_1P_{t}v_3,v_2P_1v_3,\ldots$,
$v_2P_{t}v_3$. Let the vertices in $N\cap V(v_1P_{i}v_3)$ be kept in
the queue $T_i$ according to the order in which they appear on the
path $P_i$ from $v_1$ to $v_3$, for $1\leq i \leq t$. Similarly, let
the vertices in $N\cap V(v_2P_iv_3)$ be kept in a queue $T_{i+t}$
according to the order in which they appear on the path $P_i$ from
$v_2$ to $v_3$. We may assume that
$T_{i}=\{t_1^{i},t_2^{i},\ldots,v_3\}$ for $1\leq i \leq 2t$. (The
following Figure \ref{fig2} shows the description, in which the
crosses indicate the vertices $t_j^i$.)

\begin{figure}[h,t]
\begin{center}
\input{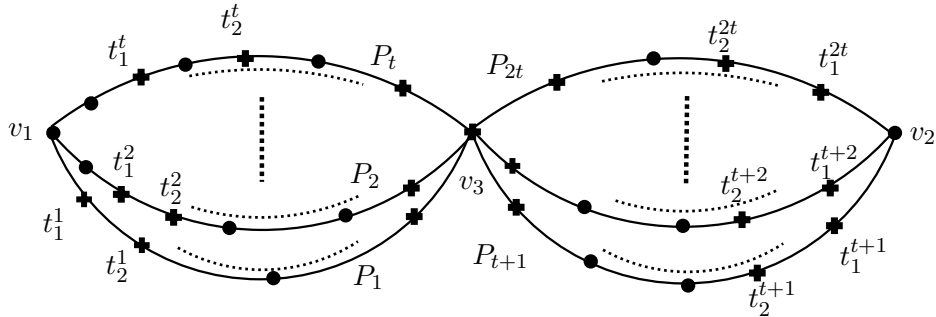}
\caption{The cross vertices of paths} \label{fig2}
\end{center}
\end{figure}

For each $t^i_j$, there exists some path $M_l$ ($1\leq l\leq 2t$)
containing $t^i_j$, since $t^i_j \in N=V(M_1\cup \cdots\cup
M_{\kappa})\cap V(P_1\cup \cdots\cup P_{t}-\{v_1,v_2\})$.

\begin{figure}[h,t]
\begin{center}
\input{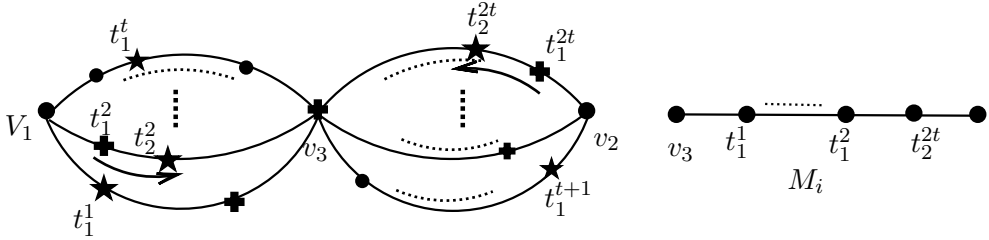}
\caption{The updated $2t$ marked vertices are $t_1^{1},
``t_2^{2}",t_1^{3},\ldots,t_1^{2t-1}, ``t_2^{2t}"$ and the updated
$2t$ marked paths are $M_{1_1},
``M_{2_2}",M_{3_1},\ldots,M_{{2t-1}_1}, ``M_{{2t}_2}"$
.}\label{fig3}
\end{center}
\end{figure}

First, we mark the vertex $t_1^{i}$ in $T_i$ for each $1\leq i \leq
2t$, and mark the corresponding path $M_{i_1}$ containing $t_1^{i}$.
If the paths $M_{1_1},M_{2_1},\ldots,M_{{2t}_1}$ are all different
(here $M_{i_j}$ denotes the corresponding path of the $j$-th vertex
of queue $T_i$), then we find $2t$ marked paths that can be used to
transform the former paths. Otherwise, there are at least two marked
vertices on the same path. Suppose that $t_1^{i_1},t_1^{i_2},\ldots$
are on the same marked path $M_i$. Keep the mark of the vertex
nearest to $v_3$ on the path $M_{i}$ and cancel the marks of the
other vertices in $\{t_1^{i_1},t_1^{i_2},\ldots\}$. Then for each
vertex $t_1^{i}$ with mark cancelled just now, mark the next vertex
$t_2^{i}$ in $T_i$ and also mark the corresponding path $M_{i_2}$
containing $t_2^{i}$.  For example, $t_1^1,t_1^2,t_1^{2t}$ are on
the same marked path, namely $M_i=M_{1_1}=M_{2_1}=M_{{2t}_1}$, and
$t_1^1$ is the vertex nearest to $v_3$ on the path $M_i$. See Figure
\ref{fig3}. (the stars indicate the updated $2t$ marked vertices.)

If the updated $2t$ marked paths are distinct, that is what we want.
Otherwise, repeat the operation like before, namely if there is a
marked path on which there are at least two marked vertices, then
cancel the marks of vertices on it except the vertex nearest to
$v_3$ and in the corresponding $T_i$ containing the vertex with mark
cancelled just now, let the next vertex and the path containing the
vertex be marked , until we find $2t$ distinct marked paths. Note
that the procedure will terminate since each $T_i$ has finite
elements and contains the special vertex $v_3$. We know that $v_3$
is a vertex of any path of $M_1,M_2,\ldots,M_{\kappa}$ (so $v_3$ can
be corresponded to any $M_i$). Therefore, if for some $T_i$, $v_3$
is marked, then we can choose anyone of the paths which has not been
marked, to be the corresponding path to mark.

There are some remarks on the procedure we described above:

\begin{rema}\label{rema1}
Finally, there are only $2t$ marked vertices $q_1,q_2,\ldots,q_{2t}$
and the $2t$ final marked vertices must be in $T_1$, $T_2$,
$\ldots$, $T_{2t}$, respectively. Since at first we choose $2t$
marked vertices which come from the $2t$ queues, respectively. Then
once we cancel the mark of a vertex, we find the next vertex to mark
in the corresponding queue. So there are always $2t$ marked vertices
which are in the $2t$ queues, respectively. Without loss of
generality, suppose $q_i\in T_i$ for $1\leq i\leq 2t$. Then
$q_i=q_j$ for $i\neq j$, if and only if $q_i=q_j=v_3$.
\end{rema}

\begin{rema}\label{rema2}
Once a path $M_i$ is marked, it will always be a marked path from
then on. Although at some step we cancel the marks of some vertices
on $M_i$, the mark of the vertex nearest to $v_3$ on $M_i$ does not
be cancelled at this step. So $M_i$ is still marked. Moreover, the
final $2t$ distinct marked paths are exactly the corresponding paths
of the $2t$ final marked vertices, respectively. Without loss of
generality, let the $2t$ distinct marked paths be
$M_{1},M_{2},\ldots,M_{2t}$ and $q_i\in V(M_{i})$.
\end{rema}

\begin{rema}\label{rema3}
If both $q_i$ and $v$ are vertices on $M_i$, $q_i$ is one of the
final marked vertices and $v$ was ever marked and then was mark
cancelled at some step, then $q_i$ is closer to $v_3$ than $v$ on
$M_i$.
\end{rema}

\begin{figure}[h,t]
\begin{center}
\input{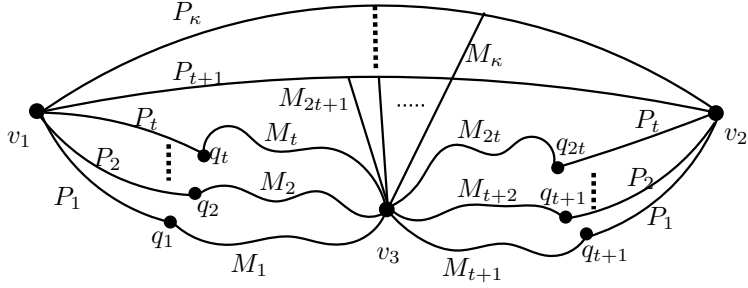}
\caption{Reduced Structure}\label{fig4}
\end{center}
\end{figure}

Now we find $2t$ marked paths $M_1,M_2,\ldots,M_{2t}$, each of which
has a final marked vertex $q_i$ such that $q_i\in T_i$, namely
$q_i,q_{i+t}\in P_i$, for $1\leq i\leq t$. Then we use the $2t$
paths to transform the former $t$ paths $P_1,\ldots,P_t$. Let
$P_{i}'=v_1P_{i}q_{i}M_{i}v_3M_{i+t}q_{i+t}P_{i}v_2$ for $1\leq
i\leq t$ (see Figure \ref{fig4}). Note that when $q_{i}=v_3$ and
$q_{i+t}=v_3$, we have $P_{i}'=P_{i}$.

\noindent {\bf Fact 1:} The $t$ walks from $v_1$ to $v_2$
$P_1',\ldots,P_t'$ are paths and have no internal vertices in common
expect $v_3$.

\pf It follows from the following three arguments:

(1) \ Since $V(\hat{v}_3M_iq_i)\subset V(M_i-v_3)$,
$V(\hat{v}_3M_jq_j)\subset V(M_j-v_3)$ and $V(M_i-v_3)\cap
V(M_j-v_3)=\emptyset$ for $i\neq j$, $V(\hat{v}_3M_iq_i)\cap
V(\hat{v}_3M_jq_j)=\emptyset$.

(2) \ Now we show that $V(\hat{q}_iP_{i-(k-1)t}\hat{v}_k)\cap
V(\hat{v}_3M_j\hat{q}_j)=\emptyset$ for $k=1$ or $2$ and $1\leq
i,j\leq 2t$. Let $v$ be a vertex in
$V(\hat{q}_{i}P_{i-(k-1)t}\hat{v}_k)$. If $v$ is on $M_j$,
obviously, $v$ is in $T_i$. Since in the queue $T_i$, $v$ is
ordered in front of $q_i$ and $q_i$ is marked, $v$ was ever
marked. So $q_j$ is closer to $v_3$ than $v$ on $M_j$ by Remark
\ref{rema2}. It follows that $v$ is not in
$V(\hat{v}_3M_j\hat{q}_j)$. If $v$ is not on $M_j$, it is
certainly not on $v_3M_jq_j$. So
$V(\hat{q}_iP_{i-(k-1)t}\hat{v}_k)\cap
V(\hat{v}_3M_j\hat{q}_j)=\emptyset$.

(3) \ It is easy to see that
$V(q_{i}P_{i-(k_1-1)t}\hat{v}_{k_1})\cap
V(q_{j}P_{j-(k_2-1)t}\hat{v}_{k_2})=\emptyset$ for $i\neq j$ and
$k_1,k_2=1$ or $2$. \qed

\noindent {\bf Fact 2:} $P_i'$ and $P_j$ are internally disjoint
paths for $1\leq i\leq t$ and $t+1\leq j\leq \kappa$.

\pf Since $P_{i}'=v_1P_{i}q_{i}M_{i}v_3M_{i+t}q_{i+t}P_{i}v_2$,
$V(v_3M_{i}q_{i})\subset V(M_i)$ and $V(v_3M_{i+t}q_{i+t})\subset
V(M_{i+t})$, obviously $V(v_3M_{i}q_{i}\cup v_3M_{i+t}q_{i+t})\cap
X=\emptyset$ for $1\leq i\leq t$. It is easy to know that
$V(\hat{v}_1P_{i}q_{i}\cup q_{i+t}P_{i}\hat{v}_2)\cap X=\emptyset$
for $1\leq i\leq t$. So $V(P_i'-\{v_1,v_2\})\cap X=\emptyset$. It
follows that $V(P_i'-\{v_1,v_2\})\cap V(P_j)=\emptyset$ for $1\leq
i\leq t$ and $t+1\leq j\leq \kappa$. \qed

\noindent {\bf Fact 3:} There is a $(\kappa-2t)$-fan from $v_3$ to
$X$ which consists of the rest paths which are not marked, namely
$\{M_{2t+1},\ldots,M_{\kappa}\}$. Moreover, the rest paths avoid the
vertices in $V(P_1'\cup\cdots \cup P_t'-\{v_1,v_2,v_3\})$

\pf By contradiction, if there exists a vertex $v$ in $V(M_i)\cap
V(P_j'-\{v_1,v_2,v_3\})$ for $2t+1\leq i\leq \kappa$ and $1\leq
j\leq t$, $v$ must be in $V(M_i)\cap
V(\hat{q}_{j+(k-1)t}P_j\hat{v}_k)$ for $k=1$ or $2$. Then we know
that $v$ was ever marked at some step and so was $M_i$. But by
Remark \ref{rema1}, if $M_i$ is marked, it will always be a marked
path from then on, a contradiction. \qed

\noindent {\bf Reduced Structure:} We have showed that by the
``Path-Transformation" we can get a structure we want, which is
called {\it Reduced Structure} (see Figure \ref{fig4}): There are
$\kappa$ $v_1v_2$-paths
$P_1',\ldots,P_{t}'$,$P_{t+1},\ldots,P_{\kappa}$ such that $v_3$ is
on $t$ paths $P_1',\ldots,P_{t}'$ of them for $1\leq t < \lceil
\frac{\kappa}{2}\rceil$, and except $v_3$ the $\kappa$ paths have no
internal vertices in common and in addition, there is a family of
$\kappa-2t$ internally disjoint $(v_3,X)$-paths
$\{M_{2t+1},\ldots,M_{\kappa}\}$ avoiding the vertices in
$V(P_1'\cup \cdots\cup P_{t}'-\{v_1,v_2,v_3\})$, where
$X=V(P_{t+1}\cup \cdots\cup P_{\kappa})$. Moreover, either the
terminal vertices of $M_{2t+1},\ldots,M_{\kappa}$ are on $\kappa-2t$
distinct paths of $P_{t+1},\ldots,P_{\kappa}$ or there are two
distinct terminal vertices on the same path. Note that $v_1$ and
$v_2$ can be regarded as vertices on any of the paths
$P_{t+1},\ldots,P_{\kappa}$.

There is still a special case we need to illustrate, namely, $|X|<
\kappa$, where $X=V(P_{t+1}\cup \cdots\cup P_{\kappa})$. Let $X'=X
\cup T$ such that $|X'|\geq \kappa$ and $v_3$ is not in $X'$, where
$T$ consists of the vertices adjacent to $v_1$ and not in $X$. So
there is a $\kappa$-fan from $v_3$ to $X'$. Then we can get $\kappa$
paths from $v_3$ to $X$ such that the terminal vertices of $|X|-1$
paths are the vertices in $X-\{v_1\}$ respectively, and all the
terminal vertices of the rest paths are $v_1$. So, by
``Path-Transformation", we can still get the Reduced Structure we
want.

{\it Proof of Theorem \ref{thm2}.}  At first, we prove that the
theorem is true for the case that $\kappa(G)=4k$, where $k$ is an
positive integer. The other cases can be verified similarly.

{\it Case 1:} $\kappa(G)=4k$ for $k\in \mathbb{N}^+$. We show that
$\kappa_3(G)\geq 3k$ by finding out $3k$ pairwise internally
disjoint trees connecting $S$, where $S$ consists of any three
vertices in $G$.

We may assume $S=\{v_1,v_2,v_3\}$. Since $G$ is
$\kappa$-connected, there are $\kappa$ pairwise internally
disjoint $v_1v_2$-paths $P_1,P_2,\ldots,P_{\kappa}$. Let $X=V(P_1
\cup \cdots\cup P_{\kappa})$.

Suppose $v_3$ is not in $X$. Obviously, $|X|\geq \kappa$ and so by
the Fan Lemma there exists a $\kappa$-fan
$\{M_1,M_2,\ldots,M_{\kappa}\}$ from $v_3$ to $X$. If the terminal
vertices $y_1,y_2,\ldots,y_{\kappa}$ of
$M_1,M_2,\ldots,M_{\kappa}$ can be regarded to be on the $\kappa$
paths $P_1,P_2,\ldots,P_{\kappa}$, respectively. Note that if the
terminal vertex is $v_1$ or $v_2$, it can be regarded as a vertex
contained in any of the paths $P_1,P_2,\ldots,P_{\kappa}$. So we
find $\kappa(G)=4k>3k$ pairwise internally disjoint trees
connecting $S$. Otherwise, there are two vertices on the same path
and without loss of generality, let $y_1, y_2\in V(P_1)$ such that
$y_1$ is closer to $v_1$ than $y_2$ on $P_1$. Then $G$ has
$\kappa$ pairwise internally disjoint $v_1v_2$-paths
$P_{1}'=v_1P_1y_1M_1v_3M_2y_2P_1v_2,P_2,\ldots,P_{\kappa}$ and
$v_3$ is on $P_1'$.

Suppose $v_3$ is in $X$, we know it must be on one of the paths
$P_1$, $P_2$, $\ldots$, $P_{\kappa}$.

Now, anyway, there exist $\kappa$ $v_1v_2$-paths
$P_1^1,P_2^1,\ldots,P_{\kappa}^1$ such that $v_3$ is on $t=1$ path
of them, say $P_1^1$, and the $\kappa$ paths have no internal
vertices in common.

For $X^1=V(P_{2}^1\cup \cdots\cup P_{\kappa}^1)$, by a $\kappa$-fan
from $v_3$ to $X^1$ and the operation ``Path-Transformation", we
adjust $P_1^1$ to $P_1^2$ and get a Reduced Structure, namely the
$\kappa$ $v_1v_2$-paths $P_1^2,P_2^1,\ldots,P_{\kappa}^1$ such that
$v_3$ is on $P_1^2$, the $\kappa$ paths have no internal vertices in
common and in addition, there is a $(\kappa-2)$-fan
$\{M_{3}^1,\ldots,M_{\kappa}^1\}$ from $v_3$ to $X^1$ avoiding the
vertices in $V(P_1^2-\{v_1,v_2,v_3\})$. Either the terminal vertices
$y_3^1,\ldots,y_{\kappa}^1$ of $M_{3}^1,\ldots,M_{\kappa}^1$ are on
$\kappa-2$ distinct paths of $P_2^1,\ldots,P_{\kappa}^1$, or there
are two distinct terminal vertices on the same path. For the former
case, we can easily find $\kappa(G)-1=4k-1>3k$ pairwise internally
disjoint trees connecting $S$. While for the latter case, We may
assume that $y_3^1, y_4^1\in V(P_2^1)$ and $y_3^1$ is closer to
$v_1$ than $y_4^1$ on $P_2^1$. Now, there are $\kappa$
$v_1v_2$-paths
$P_1^2,P_2^2=v_1P_2^1y_3^1M_3^1v_3M_4^1y_4^1P_2^1v_2,P_3^1\ldots,P_{\kappa}^1$
such that $v_3$ is on $t=2$ paths of them, say $P_1^2$ and $P_2^2$,
and the $\kappa$ paths have no internal vertices in common except
$v_3$. Then for $X^2=V(P_{3}^1\cup \cdots\cup P_{\kappa}^1)$, by a
$\kappa$-fan from $v_3$ to $X^2$ and the operation
``Path-Transformation", we adjust $P_1^2,P_2^2$ and get a Reduced
Structure. Repeat the procedure. Namely, if there are $\kappa$
$v_1v_2$-paths
$P_1^i,\ldots,P_{i}^i,P_{i+1}^{1},\ldots,P_{\kappa}^1$ for $1\leq i<
2k$ such that $v_3$ is on $t=i$ paths of them, say
$P_1^i,\ldots,P_{i}^i$, and the $\kappa$ paths have no internal
vertices in common except $v_3$, for $X^i=V(P_{i+1}^{1}\cup
\cdots\cup P_{\kappa}^1)$, by a family of $\kappa$ internally
disjoint $(v_3,X^i)$-paths and the operation ``Path-Transformation",
we adjust $P_1^i,\ldots,P_{i}^i$ to $P_1^{i+1},\ldots,P_{i}^{i+1}$
and get a Reduced Structure, namely the $\kappa$ $v_1v_2$-paths
$P_1^{i+1},\ldots,P_{i}^{i+1},P_{i+1}^{1},\ldots,P_{\kappa}^1$ such
that $v_3$ is on $t=i$ paths of them, say
$P_1^{i+1},\ldots,P_{i}^{i+1}$, the $\kappa$ paths have no internal
vertices in common except $v_3$ and in addition, there is a family
of $\kappa-2i$ internally disjoint $(v_3,X^i)$-paths
$\{M_{2i+1}^i,\ldots,M_{\kappa}^i\}$ avoiding the vertices in
$V(P_1^{i+1}\cup \cdots\cup P_{i}^{i+1}-\{v_1,v_2,v_3\})$.

Either the terminal vertices $y_{2i+1}^i,\ldots,y_{\kappa}^i$ of
$M_{2i+1}^i,\ldots,M_{\kappa}^i$ are on $\kappa-2i$ distinct paths
of $P_{i+1}^{1},\ldots,P_{\kappa}^1$, or there are two distinct
terminal vertices on the same path. For the former case, we can
find $3k$ pairwise internally disjoint trees connecting $S$ which
will be proved later and we call the case ``Middle Break". While
for the latter case, We may assume that $y_{2i+1}^i, y_{2i+2}^i\in
V(P_{i+1}^{1})$ and $y_{2i+1}^i$ is closer to $v_1$ than
$y_{2i+2}^i$ on $P_{i+1}^{1}$. Now, there are $\kappa$
$v_1v_2$-paths
$P_1^{i+1},\ldots,P_{i}^{i+1},P_{i+1}^{i+1}=v_1P_{i+1}^{1}y_{2i+1}^iM
_{2i+1}^iv_3M_{2i+2}^iy_{2i+2}^iP_{i+1}^{1}v_2,P_{i+2}^1\ldots,P_{\kappa}^1$
such that $v_3$ is on $t=i+1$ paths of them, say
$P_1^{i+1},\ldots,P_{i+1}^{i+1}$, and the $\kappa$ paths have no
internal vertices in common except $v_3$. The procedure will
terminate when either ``Middle Break" happens or $t=2k$ happens.
For the case that $t=2k$, we can also find $3k$ pairwise
internally disjoint trees connecting $S$ which will be proved
later and we call the case ``Final Break".

{\it Middle Break:} There are $\kappa$ $v_1v_2$-paths
$P_1,P_2,\ldots,P_{\kappa}$ such that $v_3$ is on $t$ paths
$P_1,\ldots,P_{t}$ of them for $1\leq t < 2k$, and except $v_3$ the
$\kappa$ paths have no internal vertices in common and in addition,
there is a family of $\kappa-2t$ internally disjoint $(v_3,X)$-paths
$\{M_{2t+1},\ldots,M_{\kappa}\}$ avoiding the vertices in $V(P_1\cup
\cdots\cup P_{t}-\{v_1,v_2,v_3\})$, where $X=V(P_{t+1}\cup
\cdots\cup P_{\kappa})$. Moreover, the terminal vertices
$y_{2t+1},\ldots,y_{\kappa}$ of $M_{2t+1},\ldots,M_{\kappa}$ are on
$\kappa-2t$ distinct paths of $P_{t+1},\ldots,P_{\kappa}$ and we may
let $y_{i}\in V(P_i)$ for $2t+1\leq i\leq \kappa$. Then we can find
pairwise internally disjoint trees connecting $S$. Let
$T_1=v_3P_1v_1P_{t+1}v_2$, $T_2=v_3P_1v_2P_{t+2}v_1$,$\ldots$,
$T_{2i-1}=v_3P_{i}v_1P_{t+2i-1}v_2$,
$T_{2i}=v_3P_{i}v_2P_{t+2i}v_1$, $\ldots$ and
$T_{t}=v_3P_{\lceil\frac{t}{2}\rceil}v_1P_{2t}v_2$ if $t$ is odd and
$T_{t}=v_3P_{\frac{t}{2}}v_2P_{2t}v_1$ if $t$ is even. Let
$T_j=P_j\cup M_j$ for $2t+1\leq j\leq \kappa$. Then let
$T_{l+t}=P_l$ for $\lceil\frac{t}{2}\rceil+1 \leq l\leq t$. So there
are
$t+(\kappa(G)-2t)+(t-\lceil\frac{t}{2}\rceil)=4k-\lceil\frac{t}{2}\rceil\geq
3k$ trees connecting $S$, since $t<2k$. Moreover, it is obvious that
the trees are pairwise internally disjoint.

{\it Final Break:} there are $\kappa$ $v_1v_2$-paths
$P_1,P_2,\ldots,P_{\kappa}$ such that $v_3$ is on $t=2k$ paths
$P_1,\ldots,P_{2k}$ of them, and the $\kappa$ paths have no
internal vertices in common except $v_3$. Let
$T_{2i-1}=v_3P_{i}v_1P_{2k+2i-1}v_2$ and
$T_{2i}=v_3P_{i}v_2P_{2k+2i}v_1$ for $1\leq i\leq k$. Then let
$T_{j+k}=P_j$ for $k+1\leq j\leq 2k$. Obviously, there are $3k$
trees connecting $S$ and they are pairwise internally disjoint.

In any case, for $\kappa(G)=4k$, we can always find $3k$ pairwise
internally disjoint trees connecting $S$, where $S$ consists of
any three vertices in $G$.

{\it Case 2:} $\kappa(G)=4k+1$ for $k\in \mathbb{N}$. It is obvious
that $\kappa_3(G)\geq 1$ when $\kappa(G)=1$. Then for $k>0$, by the
similar procedure, we can find out $3k+1$ pairwise internally
disjoint trees connecting $S$, where $S$ consists of any three
vertices in $G$. But in this case, the ``Middle Break" and ``Final
Break" have a little difference from Case 1.

{\it Middle Break:} The situation is the same as Case 1 except the
number of trees $T_j=P_j\cup M_j$. Since $\kappa(G)=4k+1$, there
are $\kappa(G)-2t=4k+1-2t$ internally disjoint $(v_3,X)$-paths
$M_{2t+1},\ldots,M_{\kappa}$ whose terminal vertices are on
$4k+1-2t$ distinct paths of $P_{t+1},\ldots,P_{\kappa}$ and so
there are $4k+1-2t$ trees $T_j=P_j\cup M_j$. Therefore, we find
$t+(4k+1-2t)+(t-\lceil\frac{t}{2}\rceil)=4k+1-\lceil\frac{t}{2}\rceil\geq
3k+1$ pairwise internally disjoint trees connecting $S$, since
$t<2k$.

{\it Final Break:} When $t=2k$, namely there are $\kappa$
$v_1v_2$-paths $P_1,P_2,\ldots,P_{\kappa}$ such that $v_3$ is on
$t=2k$ paths $P_1,\ldots,P_{2k}$ of them and the $\kappa$ paths have
no internal vertices in common except $v_3$, then we need to use the
``Path-Transformation" one more time to get a Reduced Structure.
More explicitly, for $X=V(P_{2k+1}\cup \cdots\cup P_{4k+1})$, by a
family of $\kappa$ internally disjoint $(v_3,X)$-paths and the
operation ``Path-Transformation", we get $\kappa$ $v_1v_2$-paths
$P_1',\ldots,P_{2k}'$,$P_{2k+1},\ldots$, $P_{\kappa}$ such that
$v_3$ is on $2k$ paths $P_1',\ldots,P_{2k}'$ of them, except $v_3$
the $\kappa$ paths have no internal vertices in common and in
addition, there is $\kappa-2t=1$ $(v_3,X)$-path $M_{4k+1}$ avoiding
the vertices in $V(P_1'\cup \cdots\cup P_{2k}'-\{v_1,v_2,v_3\})$. We
may assume that the terminal vertex of $M_{4k+1}$ is on $P_{4k+1}$.
Then, let $T_{2i-1}=v_3P_{i}'v_1P_{2k+2i-1}v_2$ and
$T_{2i}=v_3P_{i}'v_2P_{2k+2i}v_1$ for $1\leq i\leq k$. Let
$T_{j+k}=P_j'$ for $k+1\leq j\leq 2k$ and $T_{3k+1}$=$P_{4k+1}\cup
M_{4k+1}$. So there are $3k+1$ pairwise internally disjoint trees
connecting $S$.

{\it Case 3:} $\kappa(G)=4k+2$ for $k\in \mathbb{N}$. It is obvious
that $G$ is $(4k+1)$-connected and so by Case 2, $\kappa_3(G)\geq
3k+1$.

{\it Case 4:} $\kappa(G)=4k+3$ for $k\in \mathbb{N}$. The method is
still similar. But $1\leq t\leq 2k$ in the ``Middle Break" and
$t=2k+1$ in the ``Final Break".

From the above, if $\kappa(G)=4k+r$, we can find $3k+\lceil
\frac{r}{2}\rceil$ pairwise internally disjoint trees connecting
$S$, where $S$ consists of any three vertices in $G$, namely
$\kappa_3(G)\geq 3k+\lceil \frac{r}{2}\rceil$, for every two
integers $k$ and $r$ with $k\geq 0$ and $r\in\{0,1,2,3\}$.

Next, we will give graphs which attain the lower bound.

For $\kappa(G)=4k+2i$ with $i=0$ or $1$, we construct a graph $G$ as
follows: Let $Q=Y_1\cup Y_2$ be a vertex cut of $G$, where $Q$ is a
clique and $|Y_1|=|Y_2|=2k+i$. $G-Q$ has 2 components $C_1,C_2$.
$C_1=\{v_3\}$ and $v_3$ is adjacent to every vertex in $Q$;
$C_2=\{v_1\}\cup\{v_2\}\cup X$, $|X|=2k+i$, the induced graph of $X$
is an empty graph, every vertex in $X$ is adjacent to every vertex
in $Q\cup \{v_1,v_2\}$, $v_1$ is adjacent to every vertex in $Y_1$
and $v_2$ is adjacent to every vertex in $Y_2$. It can be checked
that $\kappa(G)=4k+2i$.

Let $S=\{v_1,v_2,v_3\}$ and let $\{T_1,T_2,\ldots,T_l\}$ be an
internally disjoint set of trees connecting $S$. For each $T_i$,
there must be a $v_1v_3$-path including a vertex in $Q$ and
$(T_i\cap Q) \cap (T_j\cap Q)=\emptyset$ for $1\leq i <j \leq l$. If
$T_i$ contains only one vertex in $Q$, then $v_3$ is a leaf of $T_i$
which means $T_i-v_3$ is still a tree connecting $v_1$ and $v_2$.
But we can see that every vertex in $Q$ is adjacent to only one of
$v_1$ and $v_2$. So $T_i-v_3$ must contain a vertex in $X$.
Therefore, there are at most $|X|$ trees in $\{T_1,T_2,\ldots,T_l\}$
containing only one vertex in $Q$ and the others contain at least
two vertices in $Q$. We can get that $l\leq
|X|+\lfloor\frac{|Q|-|X|}{2}\rfloor=2k+i+\lfloor\frac{2k+i}{2}\rfloor=3k+i$
and $\kappa_3(G)\leq \kappa(S)=l\leq 3k+i$. On the other hand,
$\kappa_3(G)\geq 3k+i$ by Theorem \ref{thm2}. It follows that
$\kappa_3(G)= 3k+i$, which means $G$ attains the lower bound.

For $\kappa(G)=4k+2i+1$ with $i=0$ or $1$, we construct a graph $G$
as follows: Let $Q=Y_1\cup Y_2 \cup \{y_0\}$ be a vertex cut of $G$,
where $Q$ is a clique and $|Y_1|=|Y_2|=2k+i$. $G-Q$ has 2 components
$C_1, \ C_2$. $C_1=\{v_3\}$ and $v_3$ is adjacent to every vertex in
$Q$; $C_2=\{v_1\}\cup\{v_2\}\cup X$, $|X|=2k+i$, the induced graph
of $X$ is an empty graph,  every vertex in $X$ is adjacent to every
vertex in $Q\cup \{v_1,v_2\}$, $v_1$ is adjacent to every vertex in
$Y_1$, $v_2$ is adjacent to every vertex in $Y_2$, and both $v_1$
and $v_2$ are adjacent to $y_0$. It can be checked similarly like
the above that $\kappa(G)=4k+2i+1$ and $\kappa_3(G)= 3k+i+1$, which
means $G$ attains the lower bound. \qed

\section{Bounds for planar graphs}

In this section we will study $\kappa_3(G)$ for planar graphs. More
precisely, we will give bounds of $\kappa_3(G)$ for planar graphs
and some graphs that attain the bounds.

First, we give the following lemma:

\begin{lem}\label{lem2}
Let $G$ be a connected graph with minimum degree $\delta$. Then
$\kappa_3(G)\leq \delta$. Especially, if there are two adjacent
vertices of degree $\delta$, then $\kappa_3(G)\leq \delta-1$.
\end{lem}
\pf We know that $\kappa(G)\leq \delta$ \cite{Bondy} and
$\kappa_3(G)\leq \kappa(G)$ by Theorem \ref{thm1}. So
$\kappa_3(G)\leq \delta$.

By contradiction, suppose that there are two adjacent vertices
$v_1$ and $v_2$ of degree $\delta$ and $\kappa_3(G)= \delta$.
Besides $v_1$ and $v_2$, we choose a vertex $v_3$ in
$V(G-\{v_1,v_2\})$ to get a set $S=\{v_1,v_2,v_3\}$. There exist
$\delta$ pairwise internally disjoint trees
$T_1,T_2,\ldots,T_{\delta}$ connecting $S$. Obviously, the
$\delta$ edges incident with $v_1$ must be contained in
$T_1,T_2,\ldots,T_{\delta}$ respectively, and so are the $\delta$
edges incident with $v_2$. Without loss of generality, we may
assume that the edge $v_1v_2$ is contained in $T_1$. But since
$T_1$ is a tree connecting $v_1,v_2$ and $v_3$, it must contain
another edge incident with $v_1$ or $v_2$, a contradiction. It
follows that $\kappa_3(G)\leq \delta-1$. \qed

By Kuratowski's Theorem \cite{Kuratowski}, a graph is planar if
and only if it contains no subdivision of $K_5$ or $K_{3,3}$. We
will use the theorem to prove the following lemma:

\begin{lem}\label{lem6}
For a connected planar graph $G$ with $\kappa_3(G)=k$, there are no
three vertices of degree $k$ in $G$, where $k\geq 3$.
\end{lem}
\pf By contradiction, let $v_1,v_2$ and $v_3$ be three vertices of
degree $k$. Because $\kappa_3(G)=k$, there exist $k$ pairwise
internally disjoint trees $T_1,T_2,\ldots,T_{k}$ connecting
$S=\{v_1,v_2,v_3\}$. Obviously, for any $i\in \{1,2,3\}$, the $k$
edges incident with $v_i$ are contained in $T_1,T_2,\ldots,T_{k}$,
respectively. Therefore, $v_1,v_2$ and $v_3$ are leaves of any tree
$T_i$, for $1\leq i\leq k$. It can be checked that, for every tree
$T_i$, there exists a vertex $t_i$ such that $T_i$ is a $3$-fan from
$t_i$ to $S$. Since $k\geq 3$, $T_1,T_2$ and $T_3$ exist. But
$T_1\cup T_2\cup T_3$ is a subdivision of $K_{3,3}$, a
contradiction. \qed

A $k$-connected graph $G$ is $minimally$ $k$-$connected$ if the
graph $G-e$ is not $k$-connected for any edge $e$, that is, if no
edge can be deleted. The following claim is an important lemma we
will use later.

\begin{lem}\label{lem5}
If $G$ is a minimally $3$-connected graph, then $\kappa_3(G-e)=2$
for any edge $e\in E(G)$.
\end{lem}
\pf For any edge $e\in E(G)$, $\kappa(G-e)=2$ and so
$\kappa_3(G-e)\leq 2$. Let $v_1,v_2$ and $v_3$ be any three vertices
in $G$.

{\it Case 1:} Two of the three vertices are connected by three
internally disjoint paths in $G-e$. Without loss of generality, we
may assume that there are three internally disjoint $v_1v_2$-paths
$P_1,P_2,P_3$ in $G-e$.

{\it Subcase 1.1:} The vertex $v_3$ is on one of the three
$v_1v_2$-paths. We may let $v_3\in V(P_1)$. Then in $G-e$ there
are two internally disjoint trees connecting $\{v_1,v_2,v_3\}$,
namely, $T_1=v_3P_1v_1P_2v_2$ and $T_2=v_3P_1v_2P_3v_1$.

{\it Subcase 1.2:} The vertex $v_3$ is not on any of the three
$v_1v_2$-paths. Let $X=V(P_1\cup P_2\cup P_3)$. Since $G-e$ is
$2$-connected, $v_3$ is not in $X$ and $|X|\geq 2$, then there
exists a $2$-fan $\{M_1,M_2\}$ from $v_3$ to $X$ by the Fan Lemma.
Let $y_1$ and $y_2$ be the two terminal vertices of $M_1$ and
$M_2$, respectively.

If $y_1$ and $y_2$ are on two of the three $v_1v_2$-paths, we may
let $y_1\in V(P_1)$ and $y_2\in V(P_2)$. Then in $G-e$ there are
two internally disjoint trees connecting $\{v_1,v_2,v_3\}$,
namely, $T_1=P_1\cup M_1$ and $T_2=P_2\cup M_2$.

If $y_1$ and $y_2$ are on the same path, we may let $y_1, y_2 \in
V(P_1)$ and let $y_1$ be closer to $v_1$ than $y_2$ on $P_1$. Then,
in $G-e$ there are two internally disjoint trees connecting
$\{v_1,v_2,v_3\}$, namely, $T_1=v_3M_1y_1P_1v_1P_2v_2$ and
$T_2=v_3M_2y_2P_1v_2P_3v_1$.

{\it Case 2:} For $v_1,v_2$ and $v_3$, any two vertices are
connected by only two internally disjoint paths in $G-e$. But we
know, in $G$, since $G$ is $3$-connected, any two vertices are
connected by three internally disjoint paths. Then, let $v_1$ and
$v_2$ be connected by three internally disjoint paths $P_1,P_2$
and $P_3$ in $G$. It is obvious that the edge $e$ is in $E(P_1\cup
P_2\cup P_3)$. We may assume $e\in E(P_3)$.

{\it Subcase 2.1:} $v_3$ is on either $P_1$ or $P_2$. Without loss
of generality, we may assume $v_3\in V(P_1)$. Let $X=V(P_2\cup
P_3)$. Since $G$ is $3$-connected, $v_3$ is not in $X$ and $|X|\geq
3$, then there exists a $3$-fan $\{M_1,M_2,M_3\}$ from $v_3$ to $X$
by the Fan Lemma. Then we know that the $3$-fan $\{M_1,M_2,M_3\}$
from $v_3$ to $X$ still exists in $G-e$, since $e\in E(G[X])$ which
means $e$ is not in $E(M_1\cup M_2\cup M_3)$. Then by the fan
$\{M_1,M_2,M_3\}$ and the operation ``Path-Transformation", we
adjust $P_1$ to $P_1'$ and get a Reduced Structure , namely three
internally disjoint $v_1v_2$-paths $P_1', P_2$ and $P_3$ such that
$e\in E(P_3)$, $v_3\in V(P_1')$ and in addition, there exists a
$(v_3,X)$-path $M_i$ avoiding the vertices in
$V(P_1'-\{v_1,v_2,v_3\})$, where $i\in \{1,2,3\}$.

\begin{figure}[h,t]
\begin{center}
\input{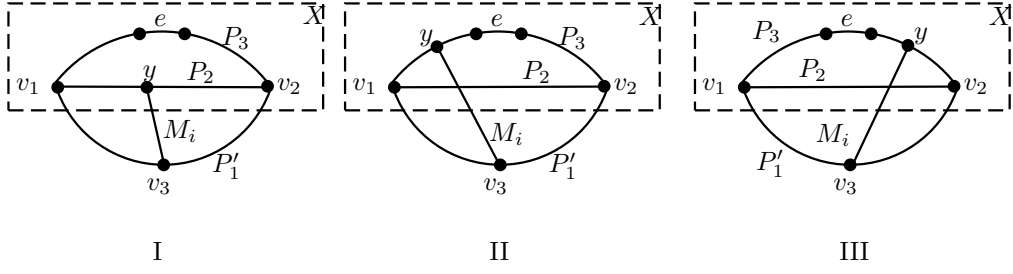}
\caption{The three graphs for Subcase 2.1} \label{fig5}
\end{center}
\end{figure}

If the terminal vertex $y$ of $M_i$ is on $P_2$, then $T_1=P_1'$ and
$T_2=P_2\cup M_i$ are two internally disjoint trees connecting
$\{v_1,v_2,v_3\}$ in $G-e$, as graph I shown in Figure \ref{fig5}.

If the terminal vertex $y$ of $M_i$ is on $P_3$ and $e\in
E(yP_3v_2)$, then $T_1=v_1P_3yM_iv_3P_1'v_2$ and
$T_2=v_3P_1'v_1P_2v_2$ are two internally disjoint trees connecting
$\{v_1,v_2,v_3\}$ in $G-e$, as graph II shown in Figure \ref{fig5}.

If the terminal vertex $y$ of $M_i$ is on $P_3$ and $e\in
E(v_1P_3y)$, then $T_1=v_2P_3yM_iv_3P_1'v_1$ and
$T_2=v_3P_1'v_2P_2v_1$ are two internally disjoint trees connecting
$\{v_1,v_2,v_3\}$ in $G-e$, as graph III shown in Figure \ref{fig5}.

{\it Subcase 2.2:} $v_3$ is on $P_3$. Without loss of generality, we
may assume $e\in E(v_3P_3v_2)$. Let $X=V(P_1\cup P_2)$. Since $G-e$
is $2$-connected, $v_3$ is not in $X$ and $|X|\geq 2$, then there
exists a $2$-fan $\{M_1,M_2\}$ from $v_3$ to $X$ by the Fan Lemma.
Let $V(M_1\cup M_2)\cap V(\hat{v}_1P_3v_3)=N$, where $N\neq
\emptyset$, since at least the vertex $v_3$ belongs to both of them.
$v$ is a vertex such that $v\in N$ and there is no vertex in $N$
closer to $v_1$ than $v$ on $P_3$, namely,
$V(\hat{v}P_3\hat{v}_1)\cap N=\emptyset$. We may let $v\in V(M_1)$
and let the terminal vertex $y$ of $M_2$ be on $P_2$. Then
$T_1=P_2\cup M_2$ and $T_2=v_3M_1vP_3v_1P_1v_2$ are two trees
connecting $\{v_1,v_2,v_3\}$ in $G-e$ and it is easy to check that
$T_1$ and $T_2$ are internally disjoint.

{\it Subcase 2.3:} $v_3$ is not on any of the three paths
$P_1,P_2,P_3$. Let $X=V(P_1\cup P_2\cup P_3)$ and then $v_3$ is
not in $X$. Since $G$ is $3$-connected and $|X|\geq 3$, there
exists a $3$-fan $\{M_1,M_2,M_3\}$ from $v_3$ to $X$ by the Fan
Lemma. We know that the $3$-fan $\{M_1,M_2,M_3\}$ from $v_3$ to
$X$ still exists in $G-e$. Let $y_1,y_2$ and $y_3$ be the terminal
vertices of $M_1,M_2$ and $M_3$, respectively.

\begin{figure}[h,t]
\begin{center}
\input{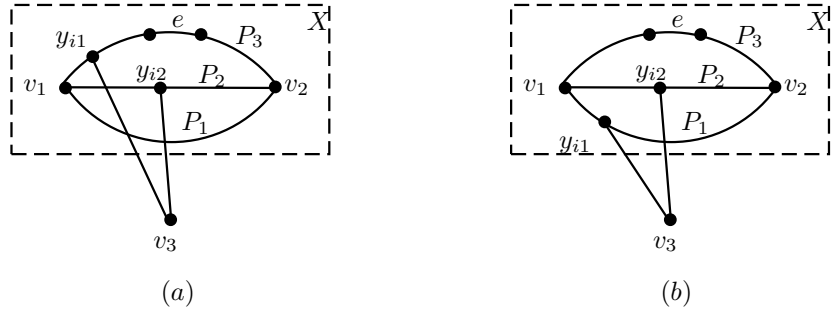}
\caption{The graphs for Subcase 2.3} \label{fig6}
\end{center}
\end{figure}

If there are two vertices $y_{i_1}$ and $y_{i_2}$ on two distinct
paths $P_{j_1}$ and $P_{j_2}$, for $1\leq i_1 \neq i_2\leq3$ and
$1\leq j_1 \neq j_2\leq3$, then it is easy to find two internally
disjoint trees connecting $\{v_1,v_2,v_3\}$ in $G-e$. See Figure
\ref{fig6}.

If the three vertices $y_1,y_2,y_3$ are on the same path $P_3$ and
$e=uv$, then either $V(v_1P_3u)$ or $V(vP_3v_2)$ contains at least
two of them. We may let $y_1$ and $y_2$ be contained in $V(v_1P_3u)$
and let $y_1$ be closer to $v_1$ than $y_2$. Then there exist three
internally disjoint $v_1v_2$-paths $P_1,P_2$ and
$P_3'=v_1P_3y_1M_1v_3M_2y_2P_3v_2$ in $G$ such that $e \in E(P_3')$
and $v_3 \in V(P_3')$, which is solved by Subcase 2.2.

If the three vertices $y_1,y_2,y_3$ are on the same path $P_1$ or
$P_2$, we may let $y_1,y_2,y_3 \in V(P_1)$ and let $y_1$ be nearest
to $v_1$ in the three vertices. Then there exist three internally
disjoint $v_1v_2$-paths $P_1'=v_1P_1y_1M_1v_3M_2y_2P_1v_2, P_2$ and
$P_3$ in $G$ such that $e \in E(P_3)$ and $v_3 \in V(P_1')$, which
is solved by Subcase 2.1.

From the above, we can always find two internally disjoint trees
connecting $\{v_1,v_2,v_3\}$ in $G-e$, where $e$ is any edge in
$G$ and $v_1,v_2,v_3$ are any three vertices in $G$. So
$\kappa_3(G-e)=2$. \qed

In the following, we list some known results which will be used
later.

\begin{lem}\label{lem8} \cite{BB}
Let $G$ be a minimally $k$-connected graph and let $T$ be the set of
vertices of degree $k$. Then $G-T$ is a (possibly empty) forest.
\end{lem}

\begin{lem}\label{lem9} \cite{BB}
A minimally $k$-connected graph of order $n$ has at least
$$\frac{(k-1)n+2}{2k-1}$$
vertices of degree $k$.
\end{lem}

\begin{lem}\label{lem7} \cite{Bondy}
Let $G$ be a $k$-connected graph and let $H$ be a graph obtained
from $G$ by adding a new vertex $y$ and joining it to at least $k$
vertices of $G$. Then $H$ is also $k$-connected.
\end{lem}

\begin{lem}\label{lem3} \cite{Bondy}
Let $G$ be a planar graph on at least three vertices. Then
$|E(G)|\leq 3|V(G)|-6$.
\end{lem}

\begin{lem}\label{lem4} \cite{Bondy}
Every planar graph has a vertex of degree at most 5, i.e.,
$\delta(G)\leq 5$.
\end{lem}

By Lemma \ref{lem4}, we only need to consider planar graphs $G$ with
connectivity $\kappa (G)$ at most 5. From Theorem \ref{thm2}, it can
be deduced that for any graph (not necessarily planar) if
$\kappa(G)=1$, $\kappa_3(G)\geq 1$; if $\kappa(G)=2$,
$\kappa_3(G)\geq 1$; if $\kappa(G)=3$, $\kappa_3(G)\geq 2$; if
$\kappa(G)=4$, $\kappa_3(G)\geq 3$, and if $\kappa(G)=5$,
$\kappa_3(G)\geq 4$. While from Theorem \ref{thm1}, we know
$\kappa_3(G)\leq \kappa(G)$, and so we get $\kappa(G)-1 \leq
\kappa_3(G)\leq \kappa(G)$, for $1\leq \kappa(G) \leq 5$. Therefore,
we get

\begin{thm}\label{thm3}
If $G$ is a connected planar graph, then $\kappa(G)-1 \leq
\kappa_3(G)\leq \kappa(G)$.
\end{thm}

Now we show that the bounds in the above theorem are sharp for
planar graphs.

{\it Case 1:} $\kappa(G)=1$. For any graph $G$ with $\kappa(G)=1$,
$\kappa_3(G)=1=\kappa(G)$. Therefore, all planar graphs with
connectivity 1 can attain the upper bound.

{\it Case 2:} $\kappa(G)=2$. There exist planar graphs with
connectivity 2 that have two adjacent vertices of degree 2. Then by
Lemma \ref{lem2}, these graphs satisfy that $\kappa_3=1$ which means
that they attain the lower bound. For example, for any cycle $C$, we
have $\kappa(C)=2$ and $\kappa_3(C)=1$.

Let $G$ be a planar minimally $3$-connected graph. By Lemma
\ref{lem5}, we know that $\kappa(G-e)=2$ and $\kappa_3(G-e)=2$ for
any edge $e\in E(G)$. Then the connected planar graph $G-e$ attains
the upper bound.

{\it Case 3:} $\kappa(G)=3$. We will show that for any planar
minimally $3$-connected graph $G$, $\kappa_3(G)=2$ which means that
it attains the lower bound.

If there are two adjacent vertices of degree 3, then by Lemma
\ref{lem2} we get $\kappa_3(G)=2$. Otherwise, any two vertices of
degree 3 are not adjacent. Let $T$ be the set of vertices of degree
3 and so $G[T]$ is an empty graph. By Lemma \ref{lem8}, we get that
$G-T$ is a forest. Let $F_1$ be a component of the forest and let
$\partial(F_1)$ denote the edge cut of $G$ associated with $V(F_1)$.
Then $|\partial(F_1)|\geq 4|F_1|-2(|F_1|-1)=2|F_1|+2>3$, since the
degree of any vertex in $V(F_1)$ is at least 4 in $G$. We know that
$N(F_1)\subseteq T$ and if there are two vertices $v_1,v_2$ in
$V(F_1)$ adjacent to a vertex $u$ in $T$ simultaneously, namely
$v_1u,v_2u \in
\partial(F_1)$, there exists a cycle $C=v_1Pv_2uv_1$, where $P$ is a
$v_1v_2$-path in $F_1$. There is just one vertex of degree 3 in
$V(C)$. But from \cite{BB} we know that each cycle of a minimally
3-connected graph contains at least two vertices of degree 3, a
contradiction. Therefore, any two vertices in $V(F_1)$ can not be
adjacent to a vertex in $T$ simultaneously, namely,
$|T|\geq|\partial(F_1)|>3$. Then in $G$ there are three vertices of
degree $3$. By Lemma \ref{lem6}, we get $\kappa_3(G)\neq 3$, namely,
$\kappa_3(G)=2$. So for any planar minimally $3$-connected graph
$G$, $\kappa_3(G)=2$ and it attains the lower bound.

Next we give graphs that attain the upper bound.

Let $G$ be a planar $4$-connected graph which is also $3$-connected
and let $G'$ be a graph obtained from $G$ by adding a new vertex $v$
to one face in some planar embedding of $G$ and joining it to $3$
vertices incident with the face. Then $G'$ is still planar and
$3$-connected by Lemma \ref{lem7}. Since there is a vertex of degree
$3$, $\kappa(G')=3$. Now we will prove $\kappa_3(G')=3$ which means
that $G'$ attains the upper bound.

For any three vertices $v_1,v_2,v_3\in V(G')$, if they are all in
$V(G)$, then restricted in $G$, we have $\kappa(\{v_1,v_2,v_3\})\geq
3$. Therefore in $G'$ it is obvious that
$\kappa(\{v_1,v_2,v_3\})\geq 3$.

Otherwise, one of them is $v$ and we may let $v_3=v$. Since
$\kappa(G)=4$, there are four internally disjoint $v_1v_2$-paths
$P_1,P_2,P_3,P_4$. Obviously, the four paths still exist in $G'$.
Let $X=V(P_1\cup P_2\cup P_3\cup P_4)$. Since $\kappa(G')=3$, $v$ is
not in $X$ and $|X|\geq 3$, then there exists a $3$-fan
$\{M_1,M_2,M_3\}$ from $v$ to $X$ by the Fan Lemma.

If the terminal vertices $y_1,y_2,y_3$ of $M_1,M_2,M_3$ are on three
of the four paths, we may let $y_1\in V(P_1), y_2\in V(P_2)$ and
$y_3\in V(P_3)$ and then there are three internally disjoint trees
connecting $\{v_1,v_2,v\}$, namely $T_1=P_1\cup M_1, T_2=P_2\cup
M_2$ and $T_3=P_3\cup M_3$.

Otherwise, there are two vertices on the same path. We may let $y_1,
y_2 \in V(P_1)$ and let $y_1$ be closer to $v_1$ than $y_2$ on
$P_1$. Then $v_1,v_2$ are connected by four internally disjoint
paths $P_1'=v_1P_1y_1M_1vM_2y_2P_1v_2, P_2,P_3,P_4$ and $v$ is on
$P_1'$. Let $X=V(P_2\cup P_3\cup P_4)$. Since $G'$ is $3$-connected,
$v$ is not in $X$ and $|X|\geq 3$, then there exists a $3$-fan
$\{M_1',M_2',M_3'\}$ from $v$ to $X$ by the Fan Lemma. Then by the
fan $\{M_1',M_2',M_3'\}$ and the operation ``Path-Transformation",
we adjust $P_1'$ to $P_1''$ and get four internally disjoint
$v_1v_2$-paths $P_1'',P_2,P_3$ and $P_4$ in $G'$ such that $v\in
V(P_1'')$ and in addition, there exists a $(v,X)$-path $M_i'$
avoiding the vertices in $V(P_1''-\{v_1,v_2,v\})$, where $i\in
\{1,2,3\}$. Without loss of generality, we may assume that the
terminal vertex $y$ of $M_i'$ is on $P_2$. Then $T_1=P_2\cup M_i',
T_2=vP_1''v_1P_3v_2$ and $T_3=vP_1''v_2P_4v_1$ are three internally
disjoint trees connecting $\{v_1,v_2,v\}$. Therefore
$\kappa_3(G')=3$. $G'$ is the graph attaining the upper bound.

{\it Case 4:} $\kappa(G)=4$. We will show that for any planar
minimally $4$-connected graph $G$, $\kappa_3(G)=3$ which means that
$G$ attains the lower bound.

Since $G$ is planar and $\kappa(G)=4$, obviously $|G|=n>5$. If
$n=6$, since $\kappa(G)=4$, the degree of any vertex is $4$ or $5$.
By Lemma \ref{lem3}, we know $|E(G)|=m\leq 3n-6$. So $4\times
6=24\leq \sum d(v)=2m\leq 6n-12=24$, which means that the degree of
every vertex is $4$. But then it is impossible that $\kappa(G)=4$.
Therefore $n\geq7$. Let $T$ be the set of vertices of degree 4.
Since $G$ is a minimally $4$-connected graph and $n\geq7$, by Lemma
\ref{lem9}, $|T|\geq \frac{3n+2}{7}>3$. Then there are three
vertices of degree $4$. By Lemma \ref{lem6}, we get $\kappa_3(G)\neq
4$, namely, $\kappa_3(G)=3$. So  any planar minimally $4$-connected
graph attains the lower bound.

It can be checked that the graphs in the following Figure \ref{fig7}
satisfy $\kappa=4$ and $\kappa_3=4$ which means that they attain the
upper bound. Moreover, we can construct a series of graphs according
to the regularity showed in Figure \ref{fig7}, which attain the
upper bound.

\begin{figure}[h,t]
\begin{center}
\input{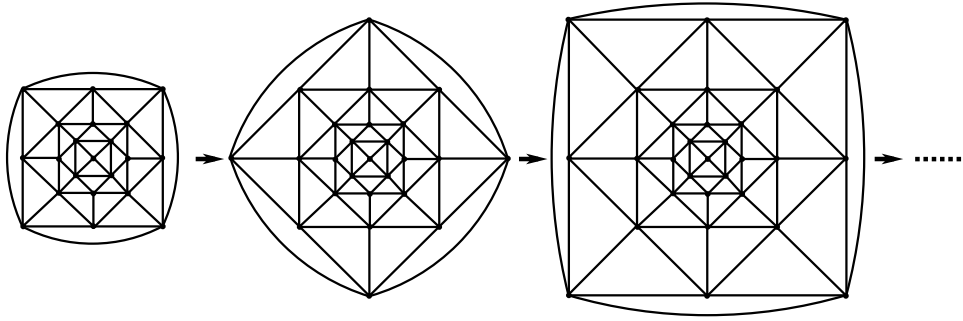}
\caption{The graphs for the upper bound of Case 4} \label{fig7}
\end{center}
\end{figure}

{\it Case 5 :} $\kappa(G)=5$. For any planar graph $G$ with
connectivity 5, if there are at most two vertices of degree $5$,
then by Lemma \ref{lem3}, $2\times 5+(n-2)\times 6\leq \sum
d(v)=2m\leq 6n-12$, namely, $6n-2\leq 6n-12$, a contradiction. So
there exist three vertices of degree $5$. By Lemma \ref{lem6}, we
get $\kappa_3(G)\neq 5$, namely, $\kappa_3(G)=4$. So, any planar
graph $G$ with connectivity 5 can attain the lower bound and
obviously can not attain the upper bound. \qed

\section{An algorithm for $\kappa_3(G)$ of planar graphs}

As well-known, for the connectivity $\kappa(G)$ of any graphs, we
have polynomial-time algorithms to get it. A natural question is
whether there is a polynomial-time algorithm to get the
$\kappa_3(G)$, or more generally, $\kappa_k(G)$. At the moment, we
do not know if such an algorithm exists for general graphs. But, for
planar graphs $G$ we shall show that $\kappa_3(G)$ can be obtained
in polynomial time, although its complexity is not very good. Since
from Theorem \ref{thm3} we have $\kappa_3(G)=\kappa (G)$ or
$\kappa(G)-1$, we only need to give a polynomial-time algorithm to
decide whether $\kappa_3(G)=\kappa (G)$.

First, it is obvious that the problem can be reduced to another
problem whether there are $\kappa(G)$ internally disjoint trees
connecting $\{v_1,v_2,v_3\}$ in polynomial time, where $v_1,v_2,v_3$
are three vertices in $V(G)$.

We now show that the problem whether there are $\kappa$ internally
disjoint trees connecting three vertices in a planar graph has a
polynomial-time algorithm.

\begin{figure}[h,t]
\begin{center}
\input{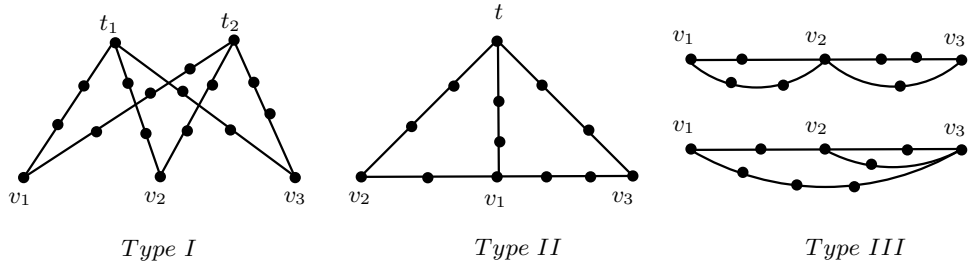}
\caption{The graphs for $T_1\cup T_2$} \label{fig8}
\end{center}
\end{figure}

For a planar graph $G$ with $\kappa(G)=2$ and three vertices
$v_1,v_2,v_3\in V(G)$, if there are two internally disjoint trees
$T_1,T_2$ connecting $\{v_1,v_2,v_3\}$, then $T_1\cup T_2$ is one of
the three types in Figure \ref{fig8}.

Our algorithm is to check all possible types until two internally
disjoint trees are found. Otherwise, we get
$\kappa(\{v_1,v_2,v_3\})=1$.

For Type I, we check for a pair of vertices $\{t_1,t_2\}\subseteq
V(G-\{v_1,v_2,v_3\})$ whether there are two internally disjoint
$3$-fans from $t_1$ to $X$ and from $t_2$ to $X$ respectively, where
$X=\{v_1,v_2,v_3\}$. If exist, we find two internally disjoint
trees. If not, we check another vertex pair until all vertex pairs
contained in $V(G-\{v_1,v_2,v_3\})$ are checked. Then we turn to
Type II.

Now the problem is that given two vertices $t_1,t_2\in
V(G-\{v_1,v_2,v_3\})$, decide whether there are two internally
disjoint $3$-fans from $t_1$ to $X$ and from $t_2$ to $X$
respectively, where $X=\{v_1,v_2,v_3\}$. At first, for each $i\in
\{1,2,3\}$, we replace the vertex $v_i$ by two new vertices
$v_{i_1}, v_{i_2}$ and let them be adjacent to all the neighbors of
$v_i$, namely, duplicating the vertex $v_i$. For each $i\in
\{1,2\}$, we replace the vertex $t_i$ by three new vertices
$t_{i_1}, t_{i_2}, t_{i_3}$ and let them be adjacent to all the
neighbors of $t_i$, namely, duplicating the vertex $t_i$ twice.
Denote the new graph by $G'$. Let $X=\{t_{1_1}, t_{1_2}, t_{1_3},
t_{2_1}, t_{2_2}, t_{2_3}\}$ and $Y=\{v_{1_1}, v_{2_1}, v_{3_1},
v_{1_2}, v_{2_2}, v_{3_2} \}$. If there exists an $XY$-linkage in
$G'$, it is easy to see that $t_{1_1}P_1v_{1_1}\cup
t_{1_2}P_2v_{2_1}\cup t_{1_3}P_3v_{3_1}$ and $t_{2_1}P_4v_{1_2}\cup
t_{2_2}P_5v_{2_2}\cup t_{2_3}P_6v_{3_2}$ can be converted into two
internally disjoint $3$-fans from $t_1$ to $\{v_1,v_2,v_3\}$ and
from $t_2$ to $\{v_1,v_2,v_3\}$ in $G$. Conversely, in $G$, any two
internally disjoint $3$-fans from $t_1$ to $\{v_1,v_2,v_3\}$ and
from $t_2$ to $\{v_1,v_2,v_3\}$ can be converted into an
$XY$-linkage in $G'$. Note that if there is an edge $e$ incident
with two vertices in $\{t_1,t_2,v_1,v_2,v_3\}$, subdivide $e$ by a
new vertex and then implement the vertex duplications. The operation
can ensure that the edge $e$ in $G$ is used only once. Since the
$k$-linkage problem, namely, the problem whether there exists an
$XY$-linkage for given sets $X$, $Y$ and any fixed value of
$|X|=|Y|=k$, has a polynomial-time algorithm, see \cite{Robertson},
then the problem whether there are two internally disjoint $3$-fans
from $t_1$ to $X$ and from $t_2$ to $X$ respectively has a
polynomial-time algorithm.

For Type II, we check for one vertex $t\in V(G-\{v_1,v_2,v_3\})$ and
the other vertex $v_{i_1}\in \{v_1,v_2,v_3\} $, whether there is a
$3$-fan from $t$ to $X$ and a $v_{i_2}v_{i_3}$-path containing
$v_{i_1}$, where the fan and the path have no vertices in common
except $v_1,v_2,v_3$, $i_1 \neq i_2 \neq i_3 \in \{1,2,3\}$ and
$X=\{v_1,v_2,v_3\}$. If exist, we find two internally disjoint trees
connecting $X$. If not, we check another vertex pair such that one
is in $V(G-\{v_1,v_2,v_3\})$ and the other is in $\{v_1,v_2,v_3\}$
until all such pairs are checked. Then we turn to Type III.

Now the problem is that given one vertex $t\in V(G-\{v_1,v_2,v_3\})$
and the other vertex $v_{i_1}\in \{v_1,v_2,v_3\}$, decide whether
there is a $3$-fan from $t$ to $X$ and a $v_{i_2}v_{i_3}$-path
containing $v_{i_1}$, where the fan and the path have no vertices in
common except $v_1,v_2,v_3$, $i_1 \neq i_2 \neq i_3 \in \{1,2,3\}$
and $X=\{v_1,v_2,v_3\}$. The method used here is the same as for
Type I. We may let $v_{i_1}=v_1$. Now for $j=2$ and $3$, replace the
vertex $v_j$ by two new vertices $v_{j_1}, v_{j_2}$ and let them be
adjacent to all the neighbors of $v_j$. Replace the vertex $t$ by
three new vertices $t_{1}, t_{2}, t_{3}$ and let them be adjacent to
all the neighbors of $t$. Replace the vertex $v_1$ by three new
vertices $v_{1_1}, v_{1_2}, v_{1_3}$ and let them be adjacent to all
the neighbors of $v_1$. Denote the new graph by $G'$. Then let
$X=\{t_1, t_2, t_3, v_{2_2}, v_{3_2}\}$ and $Y=\{v_{1_1}, v_{2_1},
v_{3_1}, v_{1_2}, v_{1_3}\}$. If there exists an $XY$-linkage in
$G'$, it is easy to see that $t_1P_1v_{1_1}\cup t_2P_2v_{2_1}\cup
t_3P_3v_{3_1}$ and $v_{2_2}P_4v_{1_2}\cup v_{1_3}P_5v_{3_2}$ can be
converted into a $3$-fan from $t$ to $\{v_1,v_2,v_3\}$ and a
$v_{2}v_{3}$-path containing $v_{1}$ in $G$ such that the fan and
the path have no vertices in common except $v_1,v_2,v_3$.
Conversely, in $G$ a $3$-fan from $t$ to $\{v_1,v_2,v_3\}$ and a
$v_{2}v_{3}$-path containing $v_{1}$ can be converted into an
$XY$-linkage in $G'$. So it can be solved in polynomial time.

For Type III, we check whether there are two $v_{i_1}v_{i_2}$-paths
both containing $v_{i_3}$, or there is a $v_{i_1}v_{i_2}$-path
containing $v_{i_3}$ and a $v_{i_1}v_{i_3}$-path containing
$v_{i_2}$. No matter what case happens, the two paths have no
vertices in common except $v_1,v_2,v_3$, $i_1 \neq i_2 \neq i_3 \in
\{1,2,3\}$ and $X=\{v_1,v_2,v_3\}$. Then the operation is similar.
Duplicate the vertices and convert the problem to the $k$-linkage
problem.

The procedure terminates when either we find $2$ internally disjoint
trees connecting $\{v_1,v_2,v_3\}$ in some type, or there are no
such two trees until all possibilities are checked. For the former
case, we get $\kappa(\{v_1,v_2,v_3\})=2$. For the latter case, we
get $\kappa(\{v_1,v_2,v_3\})=1$.

For a planar graph $G$ with $\kappa(G)=3$ and three vertices
$v_1,v_2,v_3\in V(G)$, there are at most two internally disjoint
$3$-fans from $t_1$ to $X$ and from $t_2$ to $X$ respectively, where
$t_1,t_2\in V(G-\{v_1,v_2,v_3\})$ and  $X=\{v_1,v_2,v_3\}$. So, if
$T_1, T_2, T_3$ are three internally disjoint trees connecting
$\{v_1,v_2,v_3\}$, then $T_1 \cup T_2 \cup T_3$ is one of three main
types in the following Figure \ref{fig9}.

\begin{figure}[h,t]
\begin{center}
\input{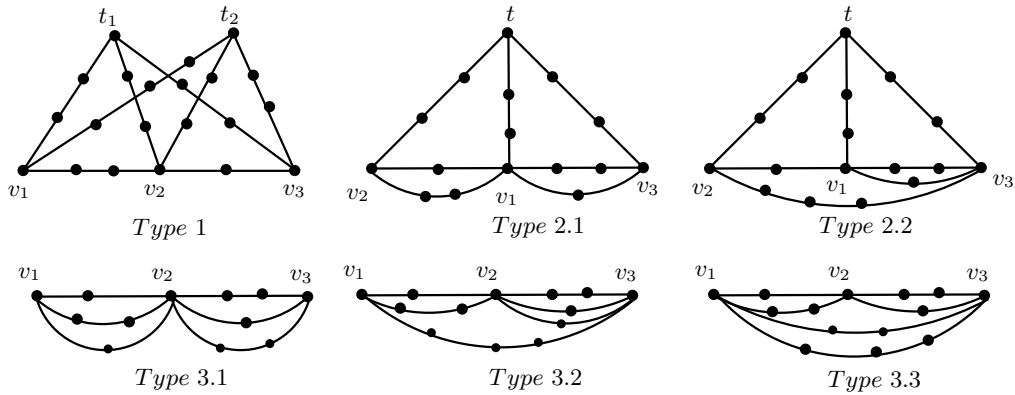}
\caption{The graphs for $T_1\cup T_2 \cup T_3$} \label{fig9}
\end{center}
\end{figure}

The method to deal with this case is similar to the case
$\kappa(G)=2$. We still implement the vertex duplications and
convert the problem into the $k$-linkage problem to solve. Still for
$\kappa(G)=4$, the method is the same. For $\kappa(G)=1$, we know
$\kappa_3(G)=\kappa(G)$ and for $\kappa(G)=5$,
$\kappa_3(G)=\kappa(G)-1$.

From the above description, we know that the algorithm is of
polynomial time. But the complexity is not very good, roughly
speaking $O(n^8)$. So how to find a more effective algorithm is an
interesting question.

\end{document}